\documentclass[12pt]{article}
\usepackage{amssymb,latexsym,amsmath}\usepackage{amsfonts}
\usepackage[cp1251]{inputenc}\usepackage[epsf]{}
\usepackage[english,russian]{babel}
\usepackage[dvips]{graphicx}
\newtheorem{thm}{Theorem}

\newtheorem{lm}{Lemma}

\begin{document}

\title{On 2-dimensional expanding attractors of A-flows\footnote{2000{\it Mathematics Subject Classification}. Primary 37D05; Secondary 37B35, 34C40}
}
\author{V.~Medvedev\and E.~Zhuzhoma}
\date{}
\maketitle

\renewcommand{\figurename}{Figure}
\renewcommand{\abstractname}{Abstract}
\renewcommand{\refname}{Bibliography}

\begin{abstract}
We prove that given any closed $n$-manifold $M^n$, $n\geq 4$, there is an A-flow $f^t$ on $M^n$ such that the non-wandering set $NW(f^t)$ consists of 2-dimensional expanding attractor (the both, orientable and non-orientable) and trivial basic sets. For 3-manifolds, we prove that given any closed 3-manifold $M^3$, there is an A-flow $f^t$ on $M^3$ such that the non-wandering set $NW(f^t)$ consists of a non-orientable 2-dimensional expanding attractor and trivial basic sets. Moreover, there is a nonsingular A-flow $f^t$ on a 3-sphere such that the non-wandering set $NW(f^t)$ consists of an orientable 2-dimensional expanding attractor and trivial basic sets (isolated periodic trajectories).
\end{abstract}

\section*{Introduction}
A-flows was introduced by Smale \cite{Smale67}. This class of flows contains structurally stable flows including Morse-Smale flows and Anosov flows.
Due to Smale's Spectral Theorem, the non-wandering set of A-flow is a disjoint union of closed transitive invariant pieces called basic sets. A basic set is called \textit{trivial} if it is either an isolated fixed point or isolated periodic trajectory. Important type of basic sets are expanding attractors introduced by Williams \cite{Williams74}. A nontrivial basic set $\Lambda$ is an \textit{expanding attractor} provided $\Lambda$ is an attractor and its topological dimension coincides with the dimension of unstable manifold at each point of $\Lambda$, see definitions at Section \ref{s:basic-definitions}. Due to Williams \cite{Williams74}, an expanding attractor consists of unstable manifolds of its points. Moreover, the unstable manifolds of points of expanding attractor form a lamination with leaves are planes and cylinders. In addition, this lamination locally homeomorphic to the product of Cantor set and Euclidean plane those dimension is more than one. Thus, the minimal topological dimension of expanding attractor for A-flows equals two, and the minimal dimension of supporting manifold equals three.
It is natural to study the following question : what manifolds admit A-flows with 2-dimensional expanding attractors ?

In the paper, we solve completely this question for closed manifolds beginning with 3-manifolds. The main results of the paper are the following statements.

\medskip
\begin{thm}\label{thm:any-M3-admitt-non-orient-2-dim}
Given any closed 3-manifold $M^3$, there is an A-flow $f^t$ on $M^3$ such that the non-wandering set $NW(f^t)$ consists of a non-orientable 2-dimensional expanding attractor and trivial basic sets.
\end{thm}
This result contrasts with the case for A-diffeomorphisms. To be precise, it follows from \cite{MedvedevZh2002,MedvedevZh2005,Plykin84,Zh82} that if a closed 3-manifold $M^3$ admits an A-diffeomorphism with 2-dimensional expanding attractor, then $\pi_1(M^3)\neq 0$. Note that there are Anosov flows with 2-dimensional expanding attractors \cite{FranksWilliams1980}. However, Margulis \cite{Margulis1967} proved that the fundamental group $\pi_1(M^3)$ of a supporting manifold $M^3$ has an exponential growth.

\medskip
\begin{thm}\label{thm:S3-admitt-orient-2-dim}
There is a nonsingular A-flow $f^t$ on a 3-sphere such that the non-wandering set $NW(f^t)$ consists of an orientable 2-dimensional expanding attractor and trivial basic sets (isolated periodic trajectories).
\end{thm}
With the help of the results above, we get the final statement.

\medskip
\begin{thm}\label{thm:any-n-manifold-admit-2-dim-attr}
Given any closed $n$-manifold $M^n$, $n\geq 4$, there is an A-flow $f^t$ on $M^n$ such that the non-wandering set $NW(f^t)$ consists of 2-dimensional expanding attractor (the both, orientable and non-orientable) and trivial basic sets.
\end{thm}

\medskip
\textit{Acknowledgments}. The authors are partially supported by Laboratory of Topological Methods in Dynamics of National Research University Higher School of Economics, of the Ministry of science and higher education of the RF, grant ag. № 075-15-2019-1931.

\section{Basic definitions}\label{s:basic-definitions}

\bigskip
\textsl{Basic definitions}.
Let $f^t$ be a smooth flow on a closed $n$-manifold $M^n$, $n\geq 3$. A subset $\Lambda\subset M^n=M$ is \textit{invariant} provided $\Lambda$ consists of trajectories of $f^t$. An invariant nonsingular set $\Lambda\subset M$ is called \textit{hyperbolic} if the sub-bundle $T_{\Lambda}M$ of the tangent bundle $TM$ can be represented as a $Df^t$-invariant continuous splitting $E^{ss}_{\Lambda}\oplus E^t_{\Lambda}\oplus E^{uu}_{\Lambda}$ such that

 1) $\dim E^{ss}_{\Lambda} + \dim E^t_{\Lambda} + \dim E^{uu}_{\Lambda} = n$;

 2) $E^t_{\Lambda}$ is the line bundle tangent to the trajectories of the flow $f^t$);

 3) there are $C_s>0$, $C_u>0$, $0<\lambda <1$ such that
$$\Vert df^t(v)\Vert \leq C_s\lambda ^t\Vert v\Vert, \quad v\in E^{ss}_{\Lambda},\quad t>0, \Vert df^{-t}(v)\Vert \leq C_u\lambda ^t\Vert v\Vert, \quad v\in E^{uu}_{\Lambda},\quad t>0.$$

If $x\in \Lambda$ is a fixed point of hyperbolic $\Lambda$, then $x$ is an isolated hyperbolic equilibrium state. The topological structure of flow near $x$ is described by Grobman-Hartman theorem, see for example \cite{Robinson-book1999}. In this case $E^t_x=0$ and $\dim E^{ss}_{\Lambda} + \dim E^{uu}_{\Lambda} = n$.

If hyperbolic $\Lambda$ does not contain fixed points, then the bundles
$$E^{uu}_{\Lambda}\oplus E^1_{\Lambda}=E^u_{\Lambda},\quad E^{ss}_{\Lambda}\oplus E^1_{\Lambda}=E^s_{\Lambda},\quad E^{uu}_{\Lambda},\quad E^{ss}_{\Lambda}$$
are uniquely integrable \cite{HirchPughShub70}, \cite{Smale67}. The corresponding leaves
$$W^u(x),\quad W^s(x),\quad W^{uu}(x),\quad W^{ss}(x)$$
through a point $x\in\lambda$ are called \textit{unstable}, \textit{stable}, \textit{strongly unstable}, and \textit{strongly stable manifolds}.

Given a set $U\subset M^n$, denote by $f^{t_0}(U)$ the shift of $A$ along the trajectories of $f^t$ on the time $t_0$. Recall that a point $x$ is non-wandering if given any neighborhood $U$ of $x$ and a number $T_0$, there is $t_0\geq T_0$ such that $U\cap f^{t_0}(U)\neq\emptyset$. The \textit{non-wandering set} $NW(f^t)$ of $f^t$ is the union of all non-wandering point.

Denote by $Fix~(f^t)$ the set of fixed points of flow $f^t$. Following Smale \cite{Smale67}, we'll call $f^t$ an \textit{A-flow} provided its non-wandering set $NW(f^t)$ is hyperbolic and the periodic trajectories are dense in
$NW(f^t) \setminus Fix~(f^t)$. Well-knows \cite{Smale67,PughShub70} that the non-wandering set $NW(f^t)$ of A-flow $f^t$ is a disjoint union of closed, and invariant, and transitive sets called \textit{basic sets}. Following Williams \cite{Williams74}, we'll call a basic set $\Omega$ an \textit{expanding attractor} provided $\Omega$ is an attractor and its topological dimension equals the dimension of unstable manifold $W^u(x)$ for every points
$x\in \Omega$. A basic set $\Lambda$ is called \textit{orientable} provided the both fiber bundles $E^{ss}_{\Lambda}$ and $E^{uu}_{\Lambda}$ are orientable. Note that if $E^{ss}_{\Lambda}$ and $E^{uu}_{\Lambda}$ are one-dimensional, then the orientability of $\Lambda$ means that the both $E^{ss}_{\Lambda}$ and $E^{uu}_{\Lambda}$ can be embedded in vector fields on $M^n$.

\section{Proof of main results}

We begin with previous results which are interesting itself. Recall that any closed manifold admits a Morse-Smale flow with a source that is a repelling fixed point. In particular, any closed manifold admits a gradient-like Morse-Smale flow having at least one source and sink \cite{Smale60a}. The following result says that any closed 3-manifold admits a Morse-Smale flow with one-dimensional repelling periodic trajectory.

\begin{lm}\label{lm:periodic-repeller-trajectory}
Given any closed 3-manifolds $M^3$, there is a Morse-Smale flow $f^t$ with repelling isolated periodic trajectory on $M^3$.
\end{lm}
\textsl{Proof}. Take a gradient-like Morse-Smale flow $f^t_0$ with a sink, say $\omega$, on $M^3$. Let $U(\omega)$ be a neighborhood of $\omega$. Without loss of generality, we can suppose that $U(\omega)$ is a ball. Since $\omega$ is a hyperbolic fixed point, one can assume that the boundary $\partial U(\omega)$ is a smooth sphere that is transversal to the trajectories. This means that the vector field inducing the flow $f^t_0$ is directed inside of $U(\omega)$. Let us introduce coordinates $(x,y,z)$ and the corresponding cylinder coordinates $(\rho,\phi,z)$ in $U(\omega)$ smoothly connected with the original coordinates in $U(\omega)$. Consider the system
$$ \dot{\rho}=\rho\cdot(1-\rho),\quad \dot{\phi}=1,\quad \dot{z}=-z. $$
It easy to check that this system has the attractive hyperbolic one-dimensional trajectory and the saddle fixed point at the origin. This complete the proof.
$\Box$

\begin{lm}\label{lm:exist-2-dim-on-s2-times-s1}
There is an A-flow on $S^2\times S^1$ such that the spectral decomposition of $f^t$ consists of two-dimensional (non-orientable) expanding attractor $\Lambda_a$ and four isolated hyperbolic repelling trajectories. Moreover, there is a neighborhood $P$ of $\Lambda_a$ homeomorphic to a solid torus $S^1\times D^2$ such that the boundary $\partial P=S^1\times S^1$ is transversal to the trajectories which enter inside of $P$ as a time parameter increases.
\end{lm}
\textsl{Proof}. Take an A-diffeomorphism $f: S^2\to S^2$ the spectral decomposition of those consists of a Plykin attractor $\Lambda_0$ and four hyperbolic sources. Due to Plykin \cite{Plykin84}, such diffeomorphism exists. Moreover, well-known that $\Lambda_0$ is an expanding attractor. Without loss of generality, one can assume that $f$ is a preserving orientation diffeomorphism (otherwise, one takes $f^2$).

Let $sus^t(f)$ be the dynamical suspension over $f$. Since $f$ is an A-diffeomorphism, $sus^t(f)$ is an A-flow. Obviously, the spectral decomposition of $f$ corresponds to the spectral decomposition of $sus^t(f)$.
Because of $f$ preserves orientation, the supporting manifold for $sus^t(f)$ is homeomorphic to $S^2\times S^1$. Since $\Lambda_0$ is a one-dimensional expanding attractor, $sus^t(f)$ has a two dimensional expanding attractor denoted by $\Lambda_a$.

Take an isolated hyperbolic repelling trajectory $\gamma$ that corresponds to some source of $f$. Since $\gamma$ is a repelling trajectory, there is a neighborhood $V(\gamma)$ homeomorphic to a solid torus such that the boundary $\partial V(\gamma)$ is transversal to the trajectories of $sus^t(f)$, so that the trajectories move outside of $V(\gamma)$ as a time parameter increases. This follows that the interior of $(S^2\times S^1)\setminus V(\gamma)$ is the neighborhood, say $P$, of $\Lambda_a$ homeomorphic to a solid torus $S^1\times D^2$ such that the boundary $\partial P=S^1\times S^1$ is transversal to the trajectories which enter inside of $P$ as a time parameter increases. Hence, $sus^t(f)=f^t$ is a desired flow.
$\Box$

\medskip
\textsc{Proof of Theorem \ref{thm:any-M3-admitt-non-orient-2-dim}}.
Let $M^3$ is a closed 3-manifold. According to Lemma \ref{lm:periodic-repeller-trajectory}, there is a Morse-Smale flow $\varphi^t$ with repelling isolated periodic trajectory $l$ on $M^3$. Hence, there is a neighborhood $U(l)$
of $l$ such that $U(l)$ is homeomorphic to the interior of the solid torus $S^1\times D^2$, and the boundary $\partial U(l)$ is transversal to the trajectories of $\varphi^t$, so that the trajectories move outside of $U(l)$ as
a time parameter increases. This follows that $U(l)$ can be replaced by the solid torus $P$ satisfying Lemma \ref{lm:exist-2-dim-on-s2-times-s1}.
As a consequence, one gets the A-flow with two-dimensional (non-orientable) expanding attractor $\Lambda_a$. This completes the proof.
$\Box$

\medskip
\textsc{Proof of Theorem \ref{thm:S3-admitt-orient-2-dim}}. The matrix
$\left(\begin{array}{cc} 2 & 1 \\ 1 & 1 \\ \end{array} \right)$
defines the Anosov diffeomorphism $A: \mathbb{T}^2\to\mathbb{T}^2$ of the two-dimensional torus $\mathbb{T}^2$. Let $f: \mathbb{T}^2\to\mathbb{T}^2$ be a DA-diffeomorphism that is constructed by modifying $A$ \cite{Robinson-book1999}. We know that the non-wandering set of $f$ consists of an orientable one-dimensional expanding attractor and isolated source. Denote by $sus^t(f)$ the dynamical suspension over $f$. Then $sus^t(f)$ is an A-flow with the non-wandering set consisting of the two-dimensional orientable expanding attractor $\Lambda$ and the isolated repelling periodic trajectory $\gamma$. The supporting manifold for $sus^t(f)$ is the mapping torus denoted by $T^3_f$. Clearly, that $T^3_f$ is a total space of a locally trivial bundle over $S^1$ and a fiber 2-torus.

Take a tubular neighborhood $U(\gamma)$ of $\gamma$ such that $U=U(\gamma)$ does not intersect $\Lambda$ and the boundary $\partial U(\gamma)=\partial U$ is transversal to the trajectories of $sus^t(f)$. Due to \cite{Christy93,Francis1983}, there is the diffeomorphism $\vartheta: \partial U\to\partial U$ such that the manifold $\left(T^3_f\setminus U\right)\bigcup_{\vartheta}U$ is a 3-sphere $S^3$. Another words, one can take out $U(\gamma)$ from $T^3_f$ and after that return $U(\gamma)$ gluing together $\left(T^3_f\setminus U\right)$ and $U(\gamma)$ along the boundary $\partial U$ under the diffeomorphism $\vartheta$. Since $\gamma$ is a repelling set, the restriction $sus^t|_{T^3_f\setminus U}$ continues to a flow $f^t$ with an isolated repelling periodic trajectory. Since $\partial U$ is transversal to $sus^t$, one can construct $f^t$ to be an A-flow. It follows that the non-wandering set $NW(f^t)$ consists of an orientable 2-dimensional expanding attractor and trivial basic set (isolated periodic trajectory). This completes the proof.
$\Box$

\medskip
\textsc{Proof of Theorem \ref{thm:any-n-manifold-admit-2-dim-attr}}. First, we construct a flow desired on an $n$-sphere $S^n$ for any $n\geq 4$. Let $f^t$ be the A-flow on $S^3$ with a 2-dimensional expanding attractor $\Lambda$. Take $S^3$ to be smoothly embedded in $S^4$ such that $S^4\setminus S^3$ is the disjoint union of 4-balls $B^4_1$, $B^4_2$. One can continue $f^t$ to $S^4$ to get A-flow, say $f^t_4$, such that $f^t_4$ has a unique source at each $B^4_i$, $i=1,2$, and $S^3$ is an attracting set for $f^t_4$. Clearly, $\Lambda$ is the 2-dimensional expanding attractor of $f^t_4$. Continuing by similar way, on can construct an A-flow, say $f^t_n$, on $S^n$ with the 2-dimensional expanding attractor $\Lambda$ for any $n\geq 5$.

Now, let $M^n$ be an arbitrary closed $n$-manifold, $n\geq 4$. Due to Smale \cite{Smale60a}, there is a gradient-like Morse-Smale flow $g^t$ on $M^n$ such that $g^t$ has a sink $s_0$. By construction, $f^t_n$ has an isolated source, say $r_n$. Take out neighborhoods $U(s_0$, $U(r_n)$ of $s_0$, $r_n$ such that the boundaries $\partial U(s_0$, $\partial U(r_n)$ are transversal to $g^t$ and $f^t_n$ respectively. If one glue $M^n\setminus U(s_0$, $S^n\setminus U(r_n)$ under a diffeomorphism $\partial U(s_0\to\partial U(r_n)$, we get a connected sum $M^n\sharp S^n$ homeomorphic to $M^n$. Then the flows $g^t$ and $f^t_n$ define the A flow desired on $M^n$.
This completes the proof.
$\Box$

\bigskip
\noindent
National Research University Higher School of Economics, 25/12 Bolshaya Pecherskaya,\\ 603005, Nizhni Novgorod, Russia

\bigskip
\noindent
\textit{E-mail:} medvedev-1942@mail.ru 

\noindent
\textit{E-mail:} zhuzhoma@mail.ru

\end{document}